\documentclass[leqno]{article}
\usepackage{amsmath}
\usepackage{amssymb}
\usepackage{amsfonts}
\usepackage{mathbbol}
\usepackage{amsthm}
\usepackage[dvips]{graphicx}
\usepackage[all]{xy}

\makeatletter
 \@addtoreset{equation}{section}
 
\makeatother

\theoremstyle{plain}
\newtheorem{thm}{Theorem}[section]

\newtheorem{lem}[thm]{Lemma}

\newtheorem{rem}[thm]{Remark}

\title{Percolation on the product graph \\ of a regular tree and a line does not satisfy 
	 \\ the triangle condition at the uniqueness threshold}
\author{Kohei Yamamoto}
\date{Tohoku university}
\begin{document}
\maketitle

	\begin{abstract}
	We consider Bernoulli bond percolation on the product graph of
	a regular tree and a line.
	We show that the triangle condition does not hold at the uniqueness threshold.
	\end{abstract}

	\section{Introduction}
	Let $G=(V,E)$ be a connected, quasi-transitive and infinite graph,
	where $V$ is the set of vertices, $E$ is the set of edges.
	In Bernoulli bond percolation, 
	each edge will be open with probability $p$,
	and closed with probability $1-p$ independently, 
	where $p \in [0,1]$ is a fixed parameter.
	Let $\Omega=\{0,1\}^E$ be the set of samples,
	where $\omega(e)=1$ means $e$ is open.
	Each $\omega \in \Omega$ is regarded as a subgraph of $G$ 
	consisting of all open edges.
	The connected components of $\omega$ are referred to as clusters.
	Let $p_c=p_c(G)$ be the critical probability for Bernoulli bond percolation 
	on $G$, that is,
	\[
	p_c= \inf \left\{ p\in[0,1] \mid 
	\text{there exists an infinite cluster almost surely} \right\},
	\]
	and let $p_u=p_u(G)$ be the uniqueness threshold for Bernoulli bond percolation 
	on $G$, that is,
	\[
	p_u= \inf \left\{ p\in[0,1] \mid 
	\text{there exists a unique infinite cluster almost surely} \right\}.
	\]
	For $p \in [0,1]$ and $x, y \in V$, let $\tau_p(x,y)$ be the probability that $x$
	and $y$ are connected in $\omega$, that is, $x$ and $y$ belong the same cluster.
	Let $\chi_p(v)$ be the expected volume of the cluster containing $v$ 
	which is defined by
	\[
	\chi_p(v) = \sum_{x \in V} \tau_p(v,x).
	\]
	This expected volume is a monotone increasing function of $p$,
	and diverges at $p_c$.
	Aizenman and Newman \cite{Aizenman} introduced the triangle condition.
	They analyzed the critical behavior of $\chi_p(v)$ if $G=\mathbb{Z}^d$ and 
	the triangle condition holds at $p_c$.
	Let $\nabla_p(v)$ be the triangle diagram which is defined by
	\[
	\nabla_p(v) = \sum_{x, y \in V} \tau_p(v,x) \tau_p(x,y) \tau_p(y,v).
	\]
	We say taht $G$ satisfies the triangle condition at $p$
	if $\nabla_p(v)<\infty$ for every $v$.
	When $G=\mathbb{Z}^d$, Hara and Slade \cite{Hara} showed that 
	the triangle condition holds at $p_c$ for all $d \ge 19$.
	This result was improved by Fitzner and van der Hofstadt \cite{Fitzner} for $d\ge11$.
	It is known $p_c=p_u$ if $G=\mathbb{Z}^d$.
	Then above result also says that the triangle condition holds at $p_u$.
	When $G$ is the product graphs,
	Kozma \cite{Kozma} showed that
	the product graph of two $d$-regular trees $T_d \Box T_d$ for $d \ge 3$
	holds the triangle condition at $p_c$.
	In 2017, Hutchcroft \cite{Hutchcroft3} showed more general cases,
	$G$ is the product graph of finitely many regular trees 
	$T_{d_1} \Box T_{d_2} \Box \cdots \Box T_{d_N}$ for $d_i \ge 3$.
	Hutchcroft \cite{Hutchcroft2} also showed that the triangle condition holds at $p_c$ 
	if $G$ is nonunimodular.
	Furthermore, Hutchcroft showed taht $p_c<p_u$ holds if $G$ is nonunimodular.
	Then we does not have the result as to whether the triangle condition holds at $p_u$.
	A nonunimodular class contains $T_d \Box \mathbb{Z}$ for $d \ge3$.
	Therefore, we consider percolation on $T_d \Box \mathbb{Z}$
	and focus on the triangle condition holds or does not.
	This graph is a vertex transitive graph.
	Then we only consider $v=o$ where $o$ is a fixed origin.
	Our main result is the following theorem.
	
	\begin{thm}
	\label{thm:triangle}
	Let $G=T_d \Box \mathbb{Z}$ for $d \ge 3$. Then we have 
	\[
	\nabla_p(o)
	\begin{cases}
	< \infty & (p<p_u)\\
	= \infty & (p=p_u).
	\end{cases}
	\] 
	\end{thm}
	To lead this result, we use a certain function $\alpha(p)$ which is defined by
	\[
	\alpha(p)=\alpha_d(p)
	=\lim_{n \to \infty} \tau_p(o, (v_n,0))^{\frac{1}{n}},
	\]
	where $v_n$ is a vertex on $T_d$ with $n$ distance from the origin.
	From a homogeneity of $T_d$, 
	$\alpha(p)$ does not depend on a choice of $v_n$.
	We abbreviate $v_n$ as $n$.
	We check on the existence of a limit.
	From FKG inequality, we have
	\[
	\tau_p(o, (n+l,0))
	\ge \tau_p(o, (n,0)) \tau_p(o, (l,0))
	\]
	for all $n,l \ge 0$.
	By using Fekete's subadditive lemma, the existence of the limit is ensured,
	and we have
	\[
	\alpha(p)
	=\lim_{n \to \infty} \tau_p(o, (n,0))^{\frac{1}{n}}
	=\sup_{n \ge 1} \tau_p(o, (n,0))^{\frac{1}{n}}.
	\]
	This function was introduced by Schonmann \cite{Schonmann},
	who showed the following inequality.
	\begin{equation}
	\label{eq:alpha(p_u)}
	\alpha(p_u) \le \frac{1}{\sqrt{b}}
	\end{equation}
	where $b=d-1$.
	By using this inequality, 
	Schonmann showed that there exists a.s. no unique infinite cluster at $p_u$.
	We will show that the equality is established, that is, 
	\begin{equation}
	\label{eq:alpha(p_u)=}
	\alpha(p_u) = \frac{1}{\sqrt{b}}.
	\end{equation}
	We introduce a example of the triangle condition, 
	let $G=T_d$, it is easy to check taht 
	\[
	\nabla_p(v)
	\begin{cases}
	< \infty & (p<\frac{1}{\sqrt{b}})\\
	= \infty & (p=\frac{1}{\sqrt{b}}).
	\end{cases}
	\] 
	We already know that $\alpha(p)$ is strictly increasing on $[0,p_u]$ 
	in \cite{Yamamoto}, that is, $\alpha(p)<1/\sqrt{b}$ for all $p<p_u$.
	We will make an upper bound and a lower bound of $\nabla_p(v)$
	by using $\alpha(p)$, and to lead Theorem \ref{thm:triangle} similar to $T_d$.

	\section{Proof}
	We define the level difference function $L(x,y)$ 
	from $T_d \times T_d$ to $\mathbb{Z}$.
	Let $\xi$ be a fixed end of $T_d$.
	The parent of a vertex $x \in T_d$ is the unique neighbor of $x$
	that is closer to $\xi$ than $x$ is. 
	We call the other vertices of $x$ its children. 
	If $y$ is parent of $x$, then we define $L(x,y)=1$.
	If $y$ is child of $x$, then we define $L(x,y)=-1$.
	In general cases, for any $x,y$,
	there exists an unique geodesic $\{x_i\}_{i=0}^n$ such that $x_0=x$ and $x_n=y$,
	then we define 
	\[
	L(x,y)= \sum_{i=1}^n L(x_{i-1},x_i).
	\]
	Note taht $L(x,z)=L(x,y)+L(y,z)$ and $L(y,x)=-L(x,y)$ for any $x,y,z \in T_d$.
	This function is extended to $T_d \Box \mathbb{Z}$ naturaly.
	Let $\pi$ be a natural projection from $T_d \Box \mathbb{Z}$ to $T_d$.
	Then we extend $L(x,y)$ as $L(x,y)=L(\pi(x), \pi(y))$.
	Similarly, we have $L(x,z)=L(x,y)+L(y,z)$ and $L(y,x)=-L(x,y)$
	for any $x,y,z \in T_d \Box \mathbb{Z}$.
	We define $\Delta(x,y)$ by 
	\[
	\Delta(x,y)=b^{L(x,y)}
	\]	
	for all $x, y \in T_d \Box \mathbb{Z}$ where $b=d-1$.
	Note that $\Delta(x,z)=\Delta(x,y)\Delta(y,z)$ and $\Delta(y,x)=\Delta(x,y)^{-1}$ 
	for any $x,y,z$.
	We define the tilted susceptibility by
	\[
	\chi_{p,1/2}(o)
	= \sum_{x \in T_d \Box \mathbb{Z}} \tau_p(o,x) \Delta(o,x)^{1/2}.
	\]
	Our method is based on \cite{Hutchcroft2}, 
	if you would like to know more detail of the tilted susceptibility,
	then please refer to \cite{Hutchcroft2}.
	Hutchcroft showed the following inequality.
	\begin{equation}
	\label{eq:tiled}
	\nabla_p(o) \le \left( \chi_{p,1/2}(o) \right)^3.
	\end{equation}
	Therefore we will show that $\chi_{p,1/2}(o)< \infty$ for $p<p_u$
	to prove the first half of Theorem \ref{thm:triangle}.
	Similar to $\alpha(p)$, the function $\beta(p)$ is defined by
	\[
	\beta(p)
	=\lim_{m \to \infty} \tau_p(o \leftrightarrow (0,m))^{\frac{1}{m}}
	=\sup_{m \ge 1} \tau_p(o \leftrightarrow (0,m))^{\frac{1}{m}}.
	\]
	By FKG inequality and the homogeneity of $T_d \Box \mathbb{Z}$, we have 
	\[
	\tau_p(o, (n,m)) \le \alpha(p)^n, \quad 
	\tau_p(o, (n,m)) \le \beta(p)^{|m|}
	\]
	for each $(n,m)$.
	For $x \in T_d$, we define $I_x(p)$ by
	\[
	I_x(p)=\sum_{m \in \mathbb{Z}} \tau_p (o, (x,m)).
	\]
	
	\begin{lem}[\cite{Yamamoto}]
	If $\alpha(p)<1/\sqrt{b}$, then we have $\beta(p)<1$.
	\end{lem}
	
	By this lemma, the function $I_x(p)$ is well-defined for $p<p_u$.
	
	\begin{lem}[\cite{Yamamoto}]
	For any $p$ such that $\alpha(p)<1/\sqrt{b}$, we have
	\[
	\lim_{|x| \to \infty} I_x(p)^{\frac{1}{|x|}} = \alpha(p).
	\]
	\end{lem}
	
	By this lemma, for any $\epsilon>0$, there exists $N \in \mathbb{N}$ such that
	\[
	I_x(p) \le (\alpha(p) + \epsilon)^{|x|}
	\]
	for any $x \in T_d$ such that $|x| \ge N$.
	For any $p<p_u$, we choose $\epsilon$ such that $\alpha(p)+\epsilon < 1/\sqrt{b}$.
	Then we have
	\begin{align*}
	\chi_{p,1/2}(o)
	&= \sum_{x \in T_d} I_x (p) \Delta(o,x)^{1/2}\\
	& \le \sum_{|x| < N} I_x (p) \Delta(o,x)^{1/2} 
	+ \sum_{|x| \ge N} (\alpha(p) + \epsilon)^{|x|} \Delta(o,x)^{1/2}.
	\end{align*}
	For $r >0$ and $z \in \mathbb{R}_{>0}$, we define the function $J(r,z)$ by
	\[
	J(r,z)= \sum_{x \in T_d} r^{|x|} z^{L(o,x)}.
	\]

	\begin{lem}[\cite{Yamamoto}]
	For any $r<1/\sqrt{b}$ and $z \in (br, 1/r)$, we have $J(r,z)<\infty$.
	\end{lem} 
	\begin{rem}
	In \cite{Yamamoto}, using level function based on the origin $o$, 
	it equal to $-L(o,x)$.
	Then $z$ appeared in \cite{Yamamoto} means $z^{-1}$ in this paper.
	\end{rem}
	By this lemma, let $r=\alpha(p)+\epsilon$ and $z=\sqrt{b}$. 
	Then we have
	\[
	 \sum_{|x| \ge N} (\alpha(p) + \epsilon)^{|x|} \Delta(o,x)^{1/2}
	\le J(r,z) <\infty.
	\]
	Therefore, we have $\nabla_p(o)<\infty$ for all $p<p_u$.
	If $p>p_u$, then there exists a constant $C(p)>0$ such that
	$\tau_p(x,y) \ge C(p)$ for all $x,y$.
	Hence, we have $\chi_{p,1/2}(o)=\infty$ for $p>p_u$.
	Hutchcroft \cite{Hutchcroft2} showed that
	the set $\left\{ p \in [0,1] \mid \chi_{p,1/2}(o)<\infty \right\}$
	is open in $[0,1]$. 
	Then we have $\chi_{p_u,1/2}(o)=\infty$.
	That means $\alpha(p)$ must equal to $1/\sqrt{b}$.
	Then we have the equation \ref{eq:alpha(p_u)=}.
	By using this result, we will show that $\nabla_{p_u}(o)=\infty$.
	Similar to $I_x(p)$, we define the function $I\hspace{-3pt}I _x(p)$ by
	\[
	I\hspace{-3pt}I _x(p) 
	= \sum_{m \in \mathbb{Z}} \tau_p(o,x)^2.
	\]
	From a homogeneity of $T_d \Box \mathbb{Z}$, 
	we have $I\hspace{-3pt}I _x(p)=I\hspace{-3pt}I _y(p)$ for any $x,y$ 
	such that $|x|=|y|$.
	For $|x|=n$, we denote $I\hspace{-3pt}I _x(p)$ as $I\hspace{-3pt}I _n(p)$.
	By using BK inequality, we have
	\[
	I\hspace{-3pt}I _{n+l}(p) 
	\le \sum_{m \in \mathbb{Z}} \sum_{k \in \mathbb{Z}}
	\tau_p(o,(n,k))^2 \tau_p(o,(l,m-k))^2
	= I\hspace{-3pt}I _n(p)I\hspace{-3pt}I _l(p)
	\]
	for all $n,l \ge 1$.
	Then by using Fekete's subadditive lemma, we have
	\[
	\inf_{n \ge 1} I\hspace{-3pt}I _n(p)^{\frac{1}{n}}
	=\lim_{n \to \infty} I\hspace{-3pt}I _n(p)^{\frac{1}{n}}.
	\]
	Since $I\hspace{-3pt}I _n(p) \ge \tau_p(o,(n,0))^2$, we have
	\[
	I\hspace{-3pt}I _n(p)^{\frac{1}{n}} \ge \tau_p(o, (n,0))^{\frac{2}{n}}.
	\]	
	By taking the limit, we have
	\[
	\lim_{n \to \infty} I\hspace{-3pt}I _n(p)^{\frac{1}{n}} \ge \alpha(p)^2
	\]
	By above equation and inequality, we obtain
	\begin{equation}
	\label{ineq:alpha^2}
	I\hspace{-3pt}I _n(p) \ge \alpha(p)^{2n}
	\end{equation}
	for all $n \ge 0$. 
	From FKG inequality, we have
	\begin{equation}
	\label{eq:square}
	\nabla_p(o)
	\ge \sum_{x,y} \tau_p(o,x)\tau_p(x,o)\tau_p(o,y)\tau_p(y,o)
	= \left( \sum_x \tau_p(o,x)^2 \right)^2.
	\end{equation}
	If $\nabla_{p_u}(o)<\infty$, then also $\sum \tau_{p_u}(o,x)^2<\infty$.
	Hence $I\hspace{-3pt}I _x(p_u)$ is well-defined.
	On the other hand, by inequality \ref{ineq:alpha^2} and equation \ref{eq:alpha(p_u)=},
	we have
	\[
	\nabla_{p_u}(o)
	= \sum_{x \in T_d} I\hspace{-3pt}I _x(p_u)
	\ge \left( \sum_{x \in T_d} \alpha(p_u)^{2|x|} \right)^2
	\ge \left( \sum_{n \ge 1} b^n \cdot \frac{1}{b^n} \right)^2
	= \infty
	\]
	Therefore, we have a contradiction,
	That means $\nabla_{p_u}(o)=\infty$.
	It ends the proof of \ref{thm:triangle}.

	{}
\noindent
{Mathematical Institute \\
 Tohoku University \\
 Sendai 980-8578 \\
 Japan\\
E-mail:kohei.yamamoto.t1@dc.tohoku.ac.jp}

	\end{document}